\input amstex
\documentstyle{amsppt}
\topmatter \magnification=\magstep1 \pagewidth{5.2 in}
\pageheight{6.7 in}
\abovedisplayskip=10pt \belowdisplayskip=10pt
\parskip=8pt
\parindent=5mm
\baselineskip=2pt
\title
 Multivariate $p$-adic $L$-function
\endtitle
\author   Taekyun Kim    \endauthor

\affil{ {\it Institute of Science Education,\\
        Kongju National University, Kongju 314-701, S. Korea\\
        e-mail: tkim64$\@$hanmail.net ( or tkim$\@$kongju.ac.kr)}}\endaffil
\thanks  2000 Mathematics Subject Classification:  11S80, 11B68, 11M99 .\endthanks
\abstract{In the recent, many mathematicians studied the multiple
zeta function in the complex number field. In this paper we
construct the $p$-adic analogue of multiple zeta function which
interpolates the generalized multiple  Bernoulli numbers attached
to $\chi$ at negative integers.
 }\endabstract
\rightheadtext{  Multivariate $p$-adic $L$-function}
\leftheadtext{T. Kim}
\endtopmatter

\document

\head \S 1. Introduction \endhead Let $p$ be a fixed prime.
Throughout this paper  $\Bbb Z_p,\,\Bbb Q_p , \,\Bbb C$ and $\Bbb
C_p$ will, respectively, denote the ring of $p$-adic rational
integers, the field of $p$-adic rational numbers, the complex
number field and the completion of algebraic closure of $\Bbb Q_p
, $ cf.[3, 10, 17].  Let $v_p$ be the normalized exponential
valuation of $\Bbb C_p$ with $|p|_p=p^{-v_p(p)}=p^{-1}.$ The
Bernoulli numbers in $\Bbb C$ are defined by
$$F(t)=\frac{t}{e^t -1}=\sum_{n=0}^{\infty}B_n \frac{t^n}{n!},
\text{ for $|t|<2\pi $}. \tag 1$$ From the definition, one has
$B_0=1 , B_1=-\frac{1}{2}, \cdots.$ Also $B_{2k+1}=0 $ for $k\geq
1$. Bernoulli numbers are used to express the special values of
Riemann zeta function $\zeta(s)=\sum_{n=1}^{\infty}\frac{1}{n^s},$
for $s\in \Bbb C ,$ namely
$\zeta(2m)=\frac{(2\pi)^{2m}(-1)^{m-1}B_2m}{2(2m)!}, $ $m\geq 1.$
For the negative integers, we note that
$\zeta(1-2m)=-\frac{B_{2m}}{2m}.$ The multiple Bernoulli numbers
of order $r$ are defined as
$$F^{r}(t)=\left(\frac{t}{e^t-1}\right)^r=\sum_{n=0}^{\infty}B_n^{(r)}\frac {t^n}{n!},
\text{ for $|t|<2\pi$, $r\in \Bbb N $}. \tag2$$
 Let $x$ be an indeterminate. Then the multiple Bernoulli
 polynomials are also defined by
$$F^{r}(t, x)=\left(\frac{t}{e^t-1}\right)^r e^{xt}=\sum_{n=0}^{\infty}B_n^{(r)}(x)\frac {t^n}{n!},
\text{ for $|t|<2\pi$, $r\in \Bbb N $}. \tag3$$ Let $\chi$ be a
primitive Dirichlet character with conductor $f\in \Bbb N$. The
generalized Bernoulli numbers attached to $\chi$, $B_{n,\chi}$,
are defined by $$F_{\chi}(t)=\sum_{a=1}^f
\frac{\chi(a)te^{at}}{e^{ft}-1}=\sum_{n=0}^{\infty}B_{n,\chi}\frac{t^n}{n!},
\text{ for $|t|<\frac{2\pi}{f}$ }. \tag4 $$ In [8, 9], the
multiple generalized Bernoulli numbers attached to $\chi$,
$B_{n,\chi}^{(r)}, $ are defined by
$$F_{\chi}^r(t)=\sum_{a_1,\cdots, a_r=1}^f\frac{\chi(\sum_{i=1}^r
a_i)t^re^{t\sum_{i=1}^r
a_i}}{(e^{ft}-1)^r}=\sum_{n=0}^{\infty}B_{n,\chi}^{(r)}\frac{t^n}{n!},
\text{ $|t|<\frac{2\pi}{f}$ }. \tag5 $$ Kubota and Leopoldt proved
the existence of meromorphic functions, $L_p(s, \chi )$, defined
over the $p$-adic number field, that serve as $p$-adic equivalents
of the Dirichlet $L$-series.
 These $p$-adic
$L$-functions interpolate the values
$$L_p(1-n, \chi)=-\frac{1}{n}(1-\chi_n(p)p^{n-1})B_{n,\chi_n},
\text{ for $n\in\Bbb N=\{1, 2,\cdots, \}  $ ,}$$ where
$B_{n,\chi}$ denote the $n$th generalized Bernoulli numbers
associated with the primitive Dirichlet character $\chi ,$ and
$\chi_n=\chi w^{-n} ,$ with $w$  the $Teichm\ddot{u}ller$
character, cf.[1-24]. In this paper we study analytic continued
function which interpolates the multiple generalized Bernoulli
numbers attached to $\chi$ at negative integers in complex plane.
Finally we will construct multivariate $p$-adic $L$-function by
using Washington's method.

 \head 2. Multivariate Dirichlet's $L$-function associated with the multiple
 generalized Bernoulli numbers attached to $\chi$  at negative integers in $\Bbb C$  \endhead

In [8, 9], the multiple generalized numbers  attached to $\chi$,
$B_{n,\chi}^{(r)}, $ are defined by
$$F_{\chi}^r(t)=\sum_{a_1,\cdots, a_r=1}^f\frac{\chi(\sum_{i=1}^r
a_i)t^re^{t\sum_{i=1}^r
a_i}}{(e^{ft}-1)^r}=\sum_{n=0}^{\infty}B_{n,\chi}^{(r)}\frac{t^n}{n!},
\text{ $|t|<\frac{2\pi}{f}$ }. $$ By (3) and (5), we easily see
that
$$B_{n,\chi}^{(r)}=f^{n-r}\sum_{a_1,\cdots,a_r=0}^{f-1}B_{n}^{(r)}(\frac{a_1+\cdots+a_r
}{f})\chi(a_1+\cdots+a_r), \text{ cf. [8, 9] }. \tag6$$ For
$s\in\Bbb C$, multiple Hurwitz's  zeta function is defined  by
$$ \zeta_r(s,x)=\sum_{n_1,\cdots,n_r=0}\frac{1}{(x+n_1 +\cdots
+n_r)^s }=\frac{1}{\Gamma (s)}\int_{0}^{\infty}F_r(-t,x)t^{s-r-1}
dt, \text{ cf. [7, 8] .} \tag7$$ Note that
$$\zeta_r(-n,x)=(-1)^r\frac{n!}{(n+r)!}B_{n+r}^{(r)}(x), \text{ for $n\in\Bbb N $, cf. [16] }. $$
We also consider the below complex integral in $\Bbb C$:
$$\frac{1}{\Gamma
(s)}\int_{0}^{\infty}F_{\chi}^r(-t)t^{s-r-1}dt=\sum_{\Sb
n_1,\cdots,n_r=0\\
n_1+\cdots+n_r\neq 0 \endSb}^{\infty}\frac{\chi(n_1+\cdots+n_r)}
{(n_1+\cdots+n_r)^s}, \text{ for $s\in\Bbb C$,
$|t|<\frac{2\pi}{f}$,} \tag8$$ where $\chi$ is the primitive
Dirichlet's character with conductor $f\in\Bbb N$.

From Eq.(8), we can derive the multivariate Dirichlet's
$L$-function in complex plane as follows:
\proclaim{Definition 1}
For $s \in \Bbb C ,$ define
$$L_r (s,\chi)=\sum_{\Sb
n_1,\cdots,n_r=0\\
n_1+\cdots+n_r\neq 0 \endSb}^{\infty}\frac{\chi(n_1+\cdots+n_r)}
{(n_1+\cdots+n_r)^s},  \tag9$$ where $\chi$ is the Dirichlet's
character with conductor $f\in\Bbb N$.\endproclaim Note that
$L_r(s,\chi)$ is meromorphic for $s\in\Bbb C$ with poles at
$s=1,2,\cdots, r.$ By (8) and (5), we obtain the following:
 \proclaim{ Theorem 2}
For $n\in\Bbb N$, we have
$$L_r (-n, \chi)=(-1)^r\frac{n!}{(n+r)!}B_{n+r,\chi}^{(r)}.\tag10$$
\endproclaim

Let $s$ be a complex variable, $a, F$ be integers with $0<a<F$.
Then we now consider the function $H_r(s;a_1, \cdots, a_r|F)$ as
follows:
$$H_r(s;a_1,\cdots,a_r|F)=\sum_{\Sb m_1,\cdots,m_r>0\\ m_i\equiv
a_i (\mod F)
\endSb}\frac{1}{(m_1+\cdots+m_r)^s}=F^{-s}\zeta_r(s,\frac{a_1+\cdots+a_r}{F}).\tag11$$
The function $H_r(s;a_1, \cdots, a_r|F)$ is a meromorphic in whole
complex plane with poles at $s=1,2, \cdots, r.$

Let $\chi (\neq 1)$ be the Dirichlet's character with conductor $F
\in \Bbb N$. Then the multivariate Dirichlet's $L$-function can be
expressed as the sum
$$L_r(s,\chi)=\sum_{a_1,\cdots,
a_r=1}^F\chi(a_1+\cdots+a_r)H_r(s;a_1,\cdots,a_r|F) \text{ for
$s\in\Bbb C$ }.\tag12$$ By simple calculation, we easily see that
$$H_r(-n;a_1,\cdots,a_r|F)=F^n(-1)^r\frac{n!}{(n+r)!}B_{n+r}^{(r)}(\frac{a_1+\cdots+a_r}{F})
\text{ for $r, n\in\Bbb N$}.\tag13$$ Thus, we have
$$L_r(-n,\chi)=(-1)^r\frac{n!}{(n+r)!}B_{n+r,\chi}^{(r)} \text{ for $n\in\Bbb N $ }. \tag14$$
By using Eq.(3), the Eq.(13) is modified  by
$$H_r(s;a_1,\cdots,a_r|F)=\frac{1}{F^r}\frac{(a_1+\cdots+a_r)^{-s+r}}{\prod_{j=1}^r(s-j)}\sum_{k=0}^{\infty}
\binom{-s+r}{k}\left(\frac{F}{a_1+\cdots+a_r} \right)^k B_k^{(r)}.
\tag15$$ From (12),(13) and (15), we can derive the below:
$$\aligned
L_r(s,\chi) &=\frac{1}{\prod_{j=1}^r
(s-j)}\frac{1}{F^r}\sum_{a_1,\cdots,a_r=1}^F\chi(a_1+\cdots
+a_r)\left(a_1+\cdots+a_r\right)^{r-s}\\
& \cdot\sum_{m=0}^{\infty}\binom{r-s}{m}
\left(\frac{a_1+\cdots+a_r}{F}\right)^m
B_m^{(r)}.\endaligned\tag16$$

Remark. The values of $L_r(s,\chi)$ at negative integers are
algebraic, hence may be regarded as lying in an extension of $\Bbb
Q_p $. We therefore look for  a $p$-adic function which agrees
with $L_r(s,\chi)$ at negative integers in the next section.

 \head 3. Multivariate $p$-adic $L$-function  \endhead

In this section we shall consider the $p$-adic  analogs of the
multivariate $L$-function, $L_r(s,\chi),$ which were introduced in
the previous section. Indeed this function is the $p$-adic
interpolation function for the generalized multiple Bernoulli
numbers attached to $\chi$. Let $w$ denote the
$Teichm\ddot{u}ller$ character, having conductor $f_w =q$. For an
arbitrary character $\chi$, we define $\chi_n=\chi w^{-n},$ where
$n\in\Bbb Z$, in the sense of the product of characters. Let
$<a>=w^{-1}(a)a=\frac{a}{w(a)}.$ Then, we note that $<a>\equiv 1$
$(\mod q\Bbb Z_p).$   Let $A_j(x)=\sum_{n=0}^{\infty}a_{n,j}x^n$,
$a_{n,j}\in\Bbb C_p$, $j=0, 1, 2,\cdots$ be a sequence of power
series, each of which converges in a fixed subset $D=\{s\in\Bbb
C_p||s|_p\leq |p^*|^{-1}p^{-\frac{1}{p-1}}\}$ of $\Bbb C_p$ such
that (1) $a_{n,j}\rightarrow a_{n, 0}$ as $j\rightarrow \infty$
for $\forall n$; (2) for each $s\in D$ and $\epsilon >0$, there
exists $n_0=n_0(s,\epsilon)$ such that $\left|\sum_{n\geq
n_0}a_{n,j}s^n\right|_p<\epsilon $ for $\forall j$. Then
$\lim_{j\rightarrow \infty}A_j(s)=A_0(s)$ for all $s\in D .$ This
is used by Washington [24] to show that each of the function
$w^{-s}(a)a^s$ and $\sum_{m=0}^{\infty}\binom sm
\left(\frac{F}{a}\right)^mB_m ,$ where $F$ is the multiple of $q$
and $f=f_{\chi}$, is analytic in $D$.

 \proclaim{ Lemma 3, cf. [24]} (Washington) Let $\chi$ be the primitive Dirichlet character, and let $F$ be a
 positive integral multiple of $q$ and $f=f_{\chi}$. Define
 $$L_p(s,\chi)=\frac{1}{s-1}\frac{1}{F}\sum_{\Sb
 a=1\\(a,p)=1\endSb}^F\chi(a)<a>^{1-s}\sum_{m=0}^{\infty}\binom{1-s}{m}\left(\frac{F}{a}\right)^mB_m
 .\tag17$$ Note that $L_p(s,\chi)$ is analytic for $s\in D$ when
 $\chi\neq 1$, and meromorphic  for $s\in D$, with simple pole at
 $s=1$ having residue $1-\frac{1}{p}$ when $\chi =1$.
Furthermore, for each $n\in\Bbb Z, $ $n\geq 1$, we have $$
L_p(1-n,\chi)=-\frac{1}{n}\left(1-\chi_n(p)p^{n-1}\right)B_{n,\chi_n}.\tag18$$
 \endproclaim
We now construct the multivariate $p$-adic $L$-function,
$L_{p,r}(s,\chi)$, which interpolates the multiple generalized
Bernoulli numbers associated with $\chi$ at negative integers.

Let $F$ be a  positive integral multiple of $q$ and $f=f_{\chi}$,
and let us define the multivariate $p$-adic $L$-function as
follows:
$$\aligned
L_{p,r}(s,\chi)=&\frac{1}{\prod_{j=1}^r (s-j)}\frac{1}{F^r}\sum_{\Sb a_1,\cdots,a_r=1\\
(a_1+\cdots+a_r,p)=1\endSb}^F\chi(a_1+\cdots+a_r)<a_1+\cdots+a_r>^{r-s}\\
&\cdot \sum_{m=0}^{\infty}\binom{r-s}{m}
\left(\frac{F}{a_1+\cdots+a_r}\right)^m
B_m^{(r)}.\endaligned\tag19$$ Then $L_{p,r}(s,\chi)$ is analytic
for $t\in\Bbb C_p$ with $|t|_p\leq 1$, provided $s\in D$, except
$s\neq 1,2,\cdots, r-1$ when $\chi \neq 1$. We now let $n\in\Bbb
Z,$ $n\geq 1$, and fix $t\in\Bbb C_p$ with $|t|_p\leq 1$. Since
$F$ must be a multiple of $f=f_{\chi_n}. $ By (6), we see that
$$B_{n+r,\chi_{n+r}}^{(r)}=F^n\sum_{a_1,\cdots,a_r=0}^{F-1}B_{n+r}^{(r)}\left(\frac{a_1+\cdots+a_r}{F}\right)
\chi_{n+r}(a_1+\cdots+a_r). \tag20$$ If $\chi_n(p)=0$, then
$(p,f_{\chi_n})=1,$ so that $\frac{F}{p}$ is multiple of
$f_{\chi_n}$. Let
$$I_0 = \{\frac{a_1+\cdots+a_r}{p}|a_1+\cdots+a_r \equiv 0 (\mod p)
\text{ for some $a_i \in \Bbb Z$ with $0\leq a_i \leq F$}\}.$$
Then we have
$$\aligned
&F^n\sum_{\Sb a_1,\cdots,a_r=0\\p|a_1+\cdots
a_r\endSb}^{F-1}B_{n+r}^{(r)}\left(\frac{a_1+\cdots+a_r}{F}\right)
\chi_{n+r}(a_1+\cdots+a_r)\\
&=p^n\left(\frac{F}{p}\right)^n\chi_n(p)\sum_{\Sb a_1,\cdots,a_r=0\\
\beta \in
I_0\endSb}^{\frac{F}{p}}\chi_{n+r}(\beta)B_{n+r}^{(r)}\left(\frac{\beta}{\frac{F}{p}}\right).
\endaligned\tag21$$
Now we also define the second multiple generalized Bernoulli
numbers attached to $\chi$ as follows:
$$B_{n+r,\chi_{n+r}}^{*(r)}= \left(\frac{F}{p}\right)^n\sum_{\Sb a_1,\cdots,a_r=0\\ \beta \in I_0\endSb}^{\frac{F}{p}}
\chi_{n+r}(\beta)B_{n+r}^{(r)}\left(\frac{\beta}{\frac{F}{p}}\right).\tag22$$
Thus, we note that
$$\aligned
&B_{n+r,\chi_{n+r}}^{(r)}-p^n\chi_n(p)B_{n+r,\chi_{n+r}}^{*(r)}\\
&=F^n\sum_{\Sb a_1,\cdots,a_r=1\\
p \nmid a_1+\cdots+a_r\endSb}^F
\chi_{n+r}(a_1+\cdots+a_r)B_{n+r}^{(r)}\left(\frac{a_1+\cdots+a_r}{F}\right).\endaligned
\tag23$$ By the definition of the multiple Bernoulli polynomials
of order $r$, we see that
$$B_{n+r}^{(r)}\left(\frac{a_1+\cdots+a_r}{F}\right)=F^{-n-r}(a_1+\cdots+a_r)^{n+r}\sum_{k=0}^n
\binom{n+r}{k}\left(\frac{F}{a_1+\cdots+a_r}\right)^k B_k^{(r)}.
\tag24$$ From (23) and (24), we can derive the below Eq.(25):
$$\aligned
B_{n+r,\chi_{n+r}}^{(r)}-p^n\chi_n(p)B_{n+r,\chi_{n+r}}^{*(r)}=&F^{-r}\sum_{\Sb a_1,\cdots, a_r=1\\
p\nmid a_1+\cdots +a_r\endSb}^F
(a_1+\cdots+a_r)^{n+r}\chi_{n+r}(a_1+\cdots+a_r)\\
&\cdot
\sum_{k=0}^{n}\binom{n+r}{k}\left(\frac{F}{a_1+\cdots+a_r}\right)^k
B_k^{(r)}. \endaligned \tag25$$ By (18) and (25), we easily see
that
$$\aligned
&L_{p,r}(-n,\chi)\\
&=\frac{(-1)^r}{\prod_{j=1}^r(n+j)}\frac{1}{F^r}\sum_{\Sb
a_1,\cdots,a_r=1\\(a_1+\cdots+a_r,p)=1\endSb}^F\chi_{n+r}(a_1+\cdots+a_r)(a_1+\cdots+a_r)^{n+r}\\
&\cdot
\sum_{m=0}^{\infty}\binom{r+m}{m}\left(\frac{F}{a_1+\cdots+a_r}
\right)^m B_m^{(r)}\\
&=(-1)^r
\frac{n!}{(n+r)!}\left(B_{n+r,\chi_{n+r}}^{(r)}-\chi_{n+r}(p)p^nB_{n+r,\chi_{n+r}}^{*(r)}\right).
\endaligned$$
Therefore we obtain the below theorem:
\proclaim{Theorem 4} Let
$F$ be a positive integral of $q$ and $f=f_{\chi}$, and let
$$\aligned
&L_{p,r}(s,\chi)\\
&=\frac{1}{\prod_{j=1}^r(s-j)}\frac{1}{F^r}\sum_{\Sb
a_1,\cdots,a_r=1\\(a_1+\cdots+a_r,p)=1\endSb}^F\chi(a_1+\cdots+a_r)<a_1+\cdots+a_r>^{r-s}\\
&\cdot\sum_{m=0}^{\infty}\binom{r-s}{m}\left(\frac{F}{a_1+\cdots+a_r}\right)^m
B_m^{(r)}.
\endaligned$$
Then, $L_{p,r}(s,\chi)$ is analytic for $t\in\Bbb C_p$, $|t|_p\leq
1 ,$ provided $s\in D$, except $s\neq 1,2, \cdots, r$. Also, if
$t\in\Bbb C_p$, $|t|_p\leq 1$, this function is analytic for $s\in
D$ when $\chi \neq 1,$ and meromorphic for $s\in D .$ Furthermore,
for each $n \in \Bbb Z$, $n\geq 1$, we have
$$L_{p,r}(-n,\chi)=(-1)^r\frac{n!}{(n+r)!}
\left(B_{n+r,\chi_{n+r}}^{(r)}-\chi_{n+r}(p)p^nB_{n+r,\chi_{n+r}}^{*(r)}\right).$$
\endproclaim
The question to construct multivariate $p$-adic $L$-function in
$p$-adic number field is still open. Theorem 4 can be considered
as a part of answer for the question which was remained in [8].

\enddocument